\magnification=\magstep1 
\baselineskip=14pt
% Reference macros. Sets up a count as a counter for refs, resets it to 0.
\newcount\refno\refno=0
% References are typed by putting \ref at the start of the line, the text
% of the reference, and a blank line afterwards.  After the blank line
% following a reference to (say) Ladner, put
%     \newcount\ladner\ladner=\refno
% Then in your document, type [\the\ladner] to obtain [XX]
% At the end of the list of references, the line :
%            \immediate\closeout\reffile
% should appear to close the references file, and then
%            \input TempReferences
% will read them in later.
\chardef\other=12
\newwrite\reffile
\immediate\openout\reffile=TempReferences
\outer\def\ref{\par\medbreak\global\advance\refno by 1
  \immediate\write\reffile{}
  \immediate\write\reffile{\noexpand\item{[\the\refno]}}
  \copytoblankline}
\def\copytoblankline{\begingroup\setupcopy\copyref}
\def\setupcopy{\def\do##1{\catcode`##1=\other}\dospecials
  \catcode`\|=\other \obeylines}
{\obeylines \gdef\copyref#1
  {\def\next{#1}%
  \ifx\next\empty\let\next=\endgroup %
  \else\immediate\write\reffile{\next} \let\next=\copyref\fi\next}}

%\ref R.J. ARCHBOLD, Topologies for primal ideals, {\sl J. London
%Math. Soc.}, (2) 36 (1987) 524-542.

%\newcount\RJA\RJA=\refno

\ref R.J. ARCHBOLD, C.J.K. BATTY, On factorial states of operator algebras, III, {\sl J. Operator Theory}, 15 (1986) 53-81.

\newcount\AB\AB=\refno

%\ref R.J. ARCHBOLD, D.W.B. SOMERSET, Quasi-standard C$^*$-algebras, {\sl Math. Proc. Camb. Phil. Soc.}, 107 (1990) 349-360.

%\newcount\AS\AS=\refno

%\ref F. BECKHOFF, Topologies on the space of ideals of a Banach algebra, {\sl Stud. Math.}, 115 (1995) 189-205.

%\newcount\Be\Be=\refno

%\ref F. BECKHOFF, Topologies of compact families on the ideal space of a Banach algebra, {\sl Stud. Math.}, 118 (1996) 63-75.

%\newcount\Bec\Bec=\refno

%\ref F. BECKHOFF, Topologies on the ideal space of a Banach algebra and spectral synthesis, {\sl Proc. Amer. Math. Soc.}, 125 (1997) 2859-2866.

%\newcount\Beck\Beck=\refno

%\ref BECKHOFF, F., An example, unpublished.

%\newcount\Beckh\Beckh=\refno

%\ref S.J. BHATT, H.V. DEDANIA, Banach algebras with unique uniform norm, {\sl %Proc. Amer. Math. Soc.}, 124 (1996) 579-584.

%\newcount\BD\BD=\refno

\ref H.G. DALES, A.M. DAVIE, Quasianalytic Banach function algebras,
{\sl J. Funct. Anal.} 13 (1973) 28-50.

\newcount\DD\DD=\refno

\ref J. DIXMIER, {\sl C$^*$-algebras}, North-Holland, Amsterdam, 1982.

\newcount\Dix\Dix=\refno

%\ref J.F. FEINSTEIN, D.W.B. SOMERSET, A note on ideal spaces of Banach algebras, {\sl Bull. London Math. Soc.}, to appear.

\ref J.F. FEINSTEIN, D.W.B. SOMERSET, Strong regularity for uniform algebras, 
to appear in {\sl Contemporary Math.}, `Proceedings of the 3rd Function Spaces Conference, Edwardsville, Illinois, May 1998', Amer. Math. Soc., Rhode Island.

\newcount\FS\FS=\refno

\ref P. GORKIN, R. MORTINI, Synthesis sets in $H^\infty + C$, preprint 
(Metz, 1998).

\newcount\GM\GM=\refno

\ref K. HOFFMAN, {\sl Banach Spaces of Analytic Functions}, Prentice-Hall, Englewood Cliffs, N.J., 1962.

\newcount\Hoff\Hoff=\refno

%\ref R. KANTROWITZ, M.M. NEUMANN, Automatic continuity of homomorphisms and derivations of algebras of continuous vector-valued functions, {\sl Czech. Math. J.}, 45(120) (1995) 747-756.

%\newcount\KN\KN=\refno

\ref J.L. KELLEY, {\sl General Topology}, Van Nostrand, Princeton, 1955.

\newcount\Kel\Kel=\refno

%\ref M. MEYER, Spectral extension property and extension of multiplicative linear functionals, {\sl Proc. Amer. Math. Soc.}, 112 (1991) 855-861.

%\newcount\Mey\Mey=\refno

\ref M.M. NEUMANN, Commutative Banach algebras and decomposable operators, {\sl Mon-atshefte Math.}, 113 (1992) 227-243.

\newcount\N\N\refno

\ref T.W. PALMER, {\sl Banach Algebras and the General Theory of $^*$-Algebras}, Vol. 1, C.U.P., New York, 1994.

\newcount\Pal\Pal=\refno

%\ref C.E. RICKART, {\sl General Theory of Banach Algebras}, Van Nostrand, London, 1960.

%\newcount\Rick\Rick=\refno

\ref W. RUDIN, Continuous functions on compact spaces without perfect subsets, {\sl Proc. Amer. Math. Soc.}, 8 (1957) 39-42.

\newcount\Ru\Ru=\refno

%\ref W. RUDIN, {\sl Functional Analysis}, Tata McGraw-Hill, New Delhi, 1973.

%\newcount\Rud\Rud=\refno

\ref W. RUDIN, {\sl Real and Complex Analysis}, Tata McGraw-Hill, New Delhi, 1974.

\newcount\Rud\Rud=\refno

\ref S. SIDNEY, More on high-order non-local uniform algebras, {\sl Illinois J. Math.}, 18 (1974) 177-192.

\newcount\Sid\Sid=\refno

\ref D.W.B. SOMERSET, Minimal primal ideals in Banach algebras, {\sl Math. Proc. Camb. Phil. Soc.}, 115 (1994) 39--52.

\newcount\MinP\MinP=\refno

%\ref D.W.B. SOMERSET, The inner derivations and the primitive ideal space of a C$^*$-algebra, {\sl J. Operator Theory}, 29 (1993) 307-321.

%\newcount\der\der=\refno

%\ref D.W.B. SOMERSET, Spectral synthesis for Banach algebras, {\sl Quart. J. Math. Oxford}, (2) 49 (1998) 501-521.

%\newcount\Syn\Syn=\refno

\ref D.W.B. SOMERSET, Ideal spaces of Banach algebras, {\sl Proc. London Math. Soc.}, (3) 78 (1999) 369-400.

\newcount\Id\Id=\refno

%\ref D.W.B. SOMERSET, Topologies on ideal spaces of Banach algebras, preprint (Aberdeen University) 1997.

%\newcount\tops\tops=\refno

\ref G. STOLZENBERG, The maximal ideal space of functions locally in an algebra, {\sl Proc. Amer. Math. Soc.}, 14 (1963) 342--345.

\newcount\Stolz\Stolz=\refno

\ref E.L. STOUT, {\sl The Theory of Uniform Algebras}, Bogden and Quigley, New York, 1971.

\newcount\St\St=\refno

\ref J. WERMER, Banach algebras and analytic functions, {\sl Advances in Math.}, 1 (1961) 51--102.
 
\newcount\We\We=\refno

\ref D. R. WILKEN, Approximate normality and function algebras on the 
interval and the circle, pages 98--111 in {\sl Function Algebras}, (Proc. Internat. Sympos. 
on Function Algebras, Tulane Univ., 1965), Scott-Foresman, Chicago, Ill., 
(1966). 
 
\newcount\Wi\Wi=\refno

\immediate\closeout\reffile

\centerline{\bf NON-REGULARITY FOR BANACH FUNCTION ALGEBRAS} 
\bigskip 
\centerline{\bf J. F. Feinstein and D. W. B. Somerset} 
\bigskip 
\bigskip 
\noindent {\bf Abstract} Let $A$ be a unital Banach function algebra with character 
space $\Phi_A$. For $x\in \Phi_A$, let $M_x$ and $J_x$ be the ideals of 
functions vanishing at $x$, and in a neighbourhood of $x$, respectively. It is 
shown that the hull of $J_x$ is connected, and that if $x$ does not belong to 
the Shilov boundary of $A$ then the set $\{ y\in\Phi_A: M_x\supseteq J_y\}$ 
has an infinite, connected subset. Various related results are given. 
\bigskip
\bigskip
\noindent {\bf Maths Reviews Classification 46J20}
\bigskip 
\bigskip 
\noindent {\bf 1. Introduction} 
\bigskip 
\noindent Let $A$ be a Banach algebra and let $Prim(A)$ be the set of 
primitive ideals of $A$. The hull-kernel topology on $Prim(A)$ is defined by 
declaring the open sets to be those of the form $\{ P\in Prim(A): 
P\not\supseteq I\}$, as $I$ varies through the closed ideals of $A$. This 
topology is compact if $A$ has an identity, but not usually Hausdorff, nor 
even $T_1$. Indeed it seems, in general, to have few useful properties, and it 
has not played a prominent part in the general theory of Banach algebras. An 
attempt to find a more useful topology has been made in [\the\Id]. 

The situation is different, however, for particular classes of Banach 
algebras, such as C$^*$-algebras and certain $L^1$-group algebras. Here the 
hull-kernel topology does have good properties such as local compactness, the 
Baire property, and (for separable C$^*$-algebras) second countability. These properties have been considerably exploited in C$^*$-algebra theory and abstract harmonic analysis. 

For commutative Banach algebras, the hull-kernel topology plays a secondary role. The primitive ideals of a (unital) commutative Banach algebra $A$ are precisely the
kernels of characters. Thus $Prim(A)$ is in bijective correspondence with the 
character space $\Phi_A$, which carries the compact, Hausdorff Gelfand 
topology. This is the topology usually employed in the study of commutative 
Banach algebras, but the hull-kernel topology (defined on $\Phi_A$ using the 
natural bijection) is also used from time to time. The hull-kernel topology is a
$T_1$ topology in this case, and is weaker than the Gelfand topology. Thus the 
two topologies coincide if and only if the hull-kernel topology is Hausdorff, 
in which case the algebra is said to be regular. Even for non-regular algebras, however, it is known that every Gelfand clopen subset of the character space is hull-kernel clopen. This is the celebrated Shilov Idempotent Theorem, see [\the\Pal; 3.5.13] for example, one of the deepest results in the whole theory. Another interesting result involving the hull-kernel topology is 
Neumann's characterization of the elements of a commutative Banach algebra 
which induce a decomposable multiplication operator---these are precisely the 
elements which are continuous with respect to the hull-kernel topology, see 
[\the\N]. 

The purpose of this paper is to analyse the failure of regularity in a non-regular algebra a little more closely. A first approach, as in C$^*$-algebra theory, might be to regard non-regularity as the failure of the Hausdorff property for the hull-kernel topology. This would lead to the study of separated points (i.e. points which can be separated by disjoint hull-kernel open sets from any point not in their closure). For a separable C$^*$-algebra, the local compactness and the second countability of the hull-kernel topology ensure the existence of a dense subset of separated points. For a separable commutative Banach algebra, however, there might not be any separated points, as the case of the disc algebra shows. Furthermore, this approach fails to take advantage of the Gelfand topology on the character space.

What we do, therefore, is to adopt the approach used in the study of spectral synthesis. Let $A$ be a Banach function algebra on a compact Hausdorff space $X$, and for $x\in X$, let $M_x$ and $J_x$ be the ideals of functions vanishing at $x$, and in a neighbourhood of $x$, respectively. The standard notion is that if $J_x$ is dense in $M_x$ then $A$ is {\sl strongly regular} at $x$. If $A$ is not strongly regular at $x$, there is still the possibility that $x$ is the only point in the hull of $J_x$, i.e. in the set $\{ y\in X:M_y\supseteq J_x\}$. In this case we will say that $x$ is an {\sl R-point}. If $x$ is not an R-point, so that the hull of $J_x$ is non-trivial, the properties of the hull of $J_x$ then become interesting. An investigation along these lines, for the non-regular algebra $H^{\infty}+C$, is conducted in [\the\GM]. One question that arose in that work, which we are able to answer, is whether the hull of $J_x$ is necessarily connected.
 
As well as the hull of $J_x$, another set which it is natural to consider is the set which we call $F_x$, defined by $F_x=\{ y\in x: M_x\supseteq J_y\}$. If $F_x$ is a singleton, we say that $x$ is a {\sl point of continuity}. In Proposition 2.2 we show that a point $x$ is a separated point in $X$ if and only if it is both an R-point and a point of continuity. Thus our general approach, given a non-regular Banach function algebra, is to ask the following questions. Firstly, how badly does regularity fail---how many R-points and non-R-points are there, and how many points of continuity and discontinuity? Secondly, if $x$ is a non-R-point, or a point of discontinuity, how large is the hull of $J_x$ or the set $F_x$? Are they finite or infinite, countable or uncountable, connected or disconnected? As we shall see, the answers to these questions can vary  depending on such things as whether the Banach function algebra is natural, or whether it is a uniform algebra.

The structure of the paper is as follows. In Section 2 we introduce the various definitions in more detail, and consider some conditions which ensure that there are an abundance of R-points, or points of continuity. In Section 3 we consider non-R-points. The main result is that if $A$ is a natural Banach function algebra then the hull of $J_x$ is connected for each point $x\in\Phi_A$. Thus the hull of $J_x$ is either a singleton, in which case $x$ is an R-point, or it is uncountable. In Section 4 we consider points of discontinuity. The main results are that if $A$ is a natural Banach function algebra and $x$ does not belong to the Shilov boundary of $A$ then $F_x$ has an infinite connected subset, while if $A$ is a uniform algebra and $x$ is a point of discontinuity then $F_x$ has a non-empty perfect subset.

Let us conclude this introduction by mentioning an interesting consequence of our work. One striking difference between C$^*$-algebras and commutative Banach algebras is that whereas there are simple examples of C$^*$-algebras with hull-kernel topology which is non-Hausdorff, but very close to being Hausdorff, the standard examples of non-regular Banach function algebras all have highly non-Hausdorff hull-kernel topologies. For example, let $A$ be the C$^*$-algebra of all sequences of two-by-two complex matrices which converge
to a diagonal matrix at infinity. Then $Prim(A)$ is isomorphic to the set of 
natural numbers, with a double point at infinity, so there are only two 
non-separated points in $Prim(A)$. The disc algebra, however, which is the most familiar non-regular Banach function algebra, has character space equal to the disc, and the hull-kernel topology on this space is only a little stronger than the cofinite topology (the weakest possible $T_1$ topology), see [\the\Hoff; p.89] for a description. Part of the motivation for this paper was the search for a non-regular Banach function algebra with a hull-kernel topology which 
was close to being Hausdorff---perhaps with only a finite or countable 
number of non-separated points. Theorem 3.2 shows, however,
that the set of points which cannot be separated from a given point is 
connected, and hence is either a singleton (the point itself) or uncountable. 
Thus if a Banach function algebra is not regular, its hull-kernel topology 
must be very far from Hausdorff. 
\bigskip 
\noindent {\bf 2. Separated points, R-points, and points of continuity} 
\bigskip 
\noindent In this section we establish various basic results about separated 
points, R-points, and points of continuity, and consider some conditions which 
ensure an abundance of such points.
\bigskip 
\noindent Let $X$ be a topological space and let $x,y\in X$. Then $x\sim y$ if 
$x$ and $y$ cannot be separated by disjoint open sets in $X$. A point of $X$ 
is a {\sl separated point} if it can be separated from every point not in its 
closure. 

Now suppose that $A$ is a Banach function algebra on a
compact, Hausdorff space $X$. In this setting, when discussing separated points, we
will always work with the hull-kernel topology on $X$. Evidently $x\sim y$ 
whenever either $M_x\supseteq J_y$ or $J_x\subseteq M_y$. 
Let us say that $x$ is an 
{\sl R-point} if for all $y \in X$ with $y\ne x$, $J_x \not\subseteq M_y$
(or in other words if, working on $X$, the hull of $J_x$ is just $\{x\}$).
We say that $x$ is a {\sl point of continuity} if for $y\ne x$, 
$J_y\not\subseteq M_x$ (or, in terms of the notation introduced
earlier, if $F_x=\{x\}$). Recall that $A$ is {\sl regular on} $X$ 
if every point of $X$ is an R-point, or 
equivalently, if every point of $X$ is a point of continuity. We say that
the algebra $A$ is {\sl regular} if is regular on $\Phi_A$.
The algebra $A$ is {\sl normal on \/} $X$ if, for every pair of disjoint 
closed sets $E$, $F$ contained in $X$, there is an $f \in A$ with $f(E) \subseteq \{ 0 \}$ and 
$f(F) \subseteq \{ 1 \}$; $A$ is {\it normal\/} if it is normal on 
$\Phi_A$. It is standard, see [\the\St; 27.2] for example, that every regular 
Banach function algebra is normal. 

With $A$ and $X$ as above, we denote the Shilov boundary of $A$ by $\Gamma_A$. 
The algebra $A$ can clearly be regarded as a
Banach function algebra on $\Phi_A$, or on $\Gamma_A$ if we wish. 
However we will also consider cases where $X$ is neither equal to the character
space nor the Shilov boundary. In the case where $X=\Phi_A$ we say that $A$ is 
{\sl natural on} $X$.
\smallskip
Unless otherwise specified, we shall only consider unital Banach function algebras.

\bigskip 
\noindent {\bf Lemma 2.1} {\sl Let $A$ be a Banach function algebra on 
a compact Hausdorff space $X$, and let $x\in X$. The following are equivalent: 

(i) $x$ is a point of continuity, 

(ii) every Gelfand neighbourhood of $x$ contains a hull-kernel neighbourhood 
of $x$, 

(iii) every net in $X$ which converges to $x$ in the hull-kernel topology 
converges to $x$ in the Gelfand topology.} 
\bigskip 
\noindent {\bf Proof.} The equivalence of (ii) and (iii) is a simple matter of 
general topology. Suppose then that (i) holds. Let $(x_{\alpha})$ be a net in 
$X$ converging to $x$ in the hull-kernel topology. Suppose for a contradiction 
that $(x_{\alpha})$ does not converge to $x$ in the Gelfand topology. Then by 
passing to a subnet if necessary, we may suppose that $(x_{\alpha})$ converges 
to $y$ in the Gelfand topology, for some $y\in X$ with $y\ne x$. Since $x$ is 
a point of continuity, there exists $f\in A$ such that $f\in J_y$ and $f(x)\ne 
0$. Hence eventually $f(x_{\alpha})=0$, since $f\in J_y$ and $(x_{\alpha})$ 
converges to $y$ in the Gelfand topology. On the other hand eventually 
$f(x_{\alpha})\ne 0$, since $(x_{\alpha})$ converges to $x$ in the hull-kernel 
topology. This contradiction shows that $(x_{\alpha})$ must converge to $x$ in 
the Gelfand topology after all. Hence (iii) holds. 

Finally, suppose that (ii) holds. Let $y\in X$ with $y\ne x$. Let $N$ be a 
Gelfand neighbourhood of $x$ such that $y$ is not in the Gelfand closure of 
$N$. By assumption there exists $f\in A$ with $f$ vanishing outside $N$ and 
such that $f(x)\ne 0$. Hence $f\in J_y$, so $M_x\not\supseteq J_y$. Thus (i) 
holds. Q.E.D. 
\bigskip 
\noindent It follows from Lemma 2.1(ii) that if $\Gamma_A \subseteq X$
then every point of continuity is contained 
in the Shilov boundary of $A$.  In particular this is the case if 
$A$ is natural on $X$, or if $A$ is a uniform algebra on $X$, 

Recall that for a Banach function algebra on a compact, Hausdorff space $X$, a 
point $x\in X$ is an {\sl independent point} if for every $\epsilon$ with 
$0<\epsilon<1$ and every neighbourhood $N$ of $x$, there exists $f\in A$ with 
$f(x)=1$ and $|f|_{X\backslash N}<\epsilon$ (where $|f|_{X\backslash N}$ is defined to be $\sup\{|f(x)|:x \in X\backslash N\}$). Lemma 2.1(ii) shows that every 
point of continuity is an independent point. 
\bigskip 
\noindent{\bf Proposition 2.2} {\sl Let $A$ be a Banach function 
algebra on a compact Hausdorff space $X$. Then a point $x\in X$ is a separated 
point if and only if it is both an R-point and a point of continuity.} 
\bigskip 
\noindent{\bf Proof.} Suppose first that $x$ is a separated point. Let $y\in 
X$ with $y\ne x$. Then there exist $f,g\in A$ with $fg=0$ and $f(x)\ne 0$, 
$g(y)\ne 0$. Thus $f\in J_y$ but $f\notin M_x$, while $g\in J_x$ but $g\notin 
M_y$. Hence $M_x\not\supseteq J_y$ and $M_y\not\supseteq J_x$. Since this is 
true for all $y\ne x$, $x$ is an R-point and a point of continuity. 

Now suppose that $x$ is both an R-point and a point of continuity. Let $y\in 
X$ with $y\ne x$. Since $x$ is an R-point, $M_y\not\supseteq J_x$, so there 
exists $g\in J_x$ such that $g(y)\ne 0$. But $x$ is a point of continuity, so 
by Lemma 2.1(ii) the Gelfand neighbourhood of $x$ on which $g$ vanishes must 
contain a hull-kernel neighbourhood of $x$. Thus there exists $f\in A$ such 
that $fg=0$, and $f(x)\ne 0$. Hence $x\not\sim y$. Since this is true for all 
$y\ne x$, $x$ is a separated point. Q.E.D. 
\bigskip 
\noindent Recall that a Banach function algebra $A$ on a compact Hausdorff 
space $X$ is {\sl weakly regular on $X$} if every non-empty 
Gelfand open subset of $X$ contains a 
non-empty hull-kernel open set. If $A$ is weakly regular and uniform on $X$, 
then $X$ is necessarily the Shilov boundary of $A$. The standard example of a 
weakly regular algebra which is not regular
is the \lq tomato-can algebra', which is the uniform 
algebra of continuous functions on a solid cylinder which are analytic on the 
base of the cylinder. This is weakly regular on its character space, which is 
the solid cylinder. 
\bigskip 
\noindent For a subset $U$ of $X$, let $\overline{U}^{h-k}$ and $\overline{U}^G$ denote the closures of $U$ in the hull-kernel and Gelfand topologies respectively.
\bigskip
\noindent {\bf Theorem 2.3} {\sl Let $A$ be a Banach function algebra on a 
compact metrizable space $X$. Then $A$ is weakly regular on $X$ if and only if 
the set of points of continuity is Gelfand dense in $X$. In this case
the set of points of continuity contains a dense $G_{\delta}$ of $X$ in the 
Gelfand topology.} 
\bigskip 
\noindent{\bf Proof.} Suppose first that set of points of continuity is  
dense in $X$ in the Gelfand topology. Then Lemma 2.1(ii) shows that $A$ is 
weakly regular on $X$. 

Conversely, suppose that $A$ is weakly regular on $X$. Let 
$(U_i)_{i=1}^{\infty}$ be a base for the Gelfand topology on $X$. For each 
$i$, the set $\overline{U_i}^{h-k}\backslash U_i$ is a hull-kernel closed set 
with no hull-kernel interior, and hence no Gelfand interior. 
Thus $Y:=\bigcup_{i=1}^{\infty} \overline{U_i}^{h-k}\backslash U_i$ is a 
Gelfand meagre subset of $X$. Let $y\in X\backslash Y$. Then for any 
$x\in X$ with $x\ne y$, 
there exists a $U_i$ containing $x$ but not $y$. Since 
$y\notin\overline{U_i}^{h-k}$, there exists a function $f\in A$ which is 
non-zero at $y$ but 
vanishes on the Gelfand neighbourhood $U_i$ of $x$. Hence $M_y\not\supseteq 
J_x$. This shows that $y$ is a point of continuity. Q.E.D. 
\bigskip 
\noindent It follows from Theorem 2.3 that if $A$ is weakly regular on a 
metrizable space $X$ then the hull of $J_x$ has empty interior for each $x\in 
X$. It would be interesting to know whether Theorem 2.3 can be improved to 
show that $X$ has to have a dense $G_{\delta}$ of separated 
points.
\bigskip
\noindent {\bf Theorem 2.4} {\sl Let $A$ be a Banach function algebra 
on a compact Hausdorff space $X$, and let $y\in X$. If the hull-kernel 
topology is first countable at $y$ then $F_y$ has no Gelfand interior. If the 
hull-kernel topology is second countable on $X$ then the set of R-points contains a dense $G_{\delta}$ of $X$ in the Gelfand topology.} 
\bigskip 
\noindent {\bf Proof.} Let $V$ be a Gelfand open subset of $\Phi_A$ not 
containing $y$. Let $(U_i)_{i=1}^{\infty}$ be a base for the hull-kernel 
topology at $y$. For each $i$, set $V_i=V\cap U_i$. Then for each $x\in V$ 
there is an $i$ such that $x\notin V_i$. Hence $\bigcap_{i=1}^{\infty} V_i$ is 
empty. Since $V$ is a Baire space, it follows that at least one $V_i$ is not 
dense in $V$ in the Gelfand topology. Thus there exists $f\in A$ such that 
$f(y)\ne 0$ but $f$ vanishes on a Gelfand open subset, $W$ say, of $V$. Thus 
for $x\in W$, $x\notin F_y$. It follows that $F_y$ has empty Gelfand interior. 

Now let $(U_i)_{i=1}^{\infty}$ be a base for the hull-kernel topology on $X$. 
Then for each $i$, $\overline{U_i}^G\backslash U_i$ is a Gelfand closed set 
with no Gelfand interior. Hence $Y:=\bigcup_{i=1}^{\infty} 
\overline{U_i}^G\backslash U_i$ is a meagre subset of $X$. Let $y\in 
X\backslash Y$. Then since $(U_i)_{i=1}^{\infty}$ is a base for the 
hull-kernel topology on $X$, for any $x\in X$ there exists a $U_i$ containing 
$x$ but not $y$. Since $y\notin\overline{ U_i}^G$, there exists a function 
$f\in A$ which is non-zero at $x$ but vanishes in a Gelfand neighbourhood of 
$y$. Hence $M_x\not\supseteq J_y$. This shows that $y$ is an R-point. Q.E.D. 
\bigskip 
\noindent Lemma 2.1(ii) shows that the hull-kernel topology is first countable 
at every point of continuity, provided that the space $X$ is first countable 
in the Gelfand topology. 
\bigskip
\noindent{\bf Example} Let $X = \{ 0, 1, {1\over 2}, {1\over 3}, \dots\}$, and let $A$ be 
the restriction to $X$ of the disc algebra. By the identity principle, the restriction map is an isomorphism, and so $A$ is a Banach function algebra on $X$ (where the norm is the uniform norm of the functions on the closed unit disc). The hull-kernel topology on $X$ is simply the cofinite topology, which is second countable because $X$ is countable. By Theorem 2.4, therefore, the set of R-points is dense in $X$. In fact every point except $0$ is an R-point, while $0$ is the only point of continuity.
\bigskip
\noindent We do not have an example of a natural, non-regular Banach function algebra for which the hull-kernel topology is second countable, or even first countable. It could well be that such things do not exist. A partial
result in this direction is given in Section 4.
\bigskip 
\bigskip 
\noindent {\bf 3. Non-R-points} 
\bigskip 
\noindent In this section we consider what happens when there is a non-R-point 
$x$. We show that the hull of $J_x$ is connected, provided one is working on 
the character space. Thus there must be uncountably many points of 
discontinuity, and also uncountably many points which cannot be separated from 
$x$ in the hull-kernel topology by disjoint open sets. We also show that, although there might not be a net in $X$ converging to each point of the hull of $J_x$ in the hull-kernel topology, if $A$ is a separable Banach function algebra then the set of points $x$ which do have this property contains a dense $G_{\delta}$ subset of $X$.
\bigskip 
\noindent {\bf Proposition 3.1} {\sl Let $A$ be a natural, unital Banach function algebra. Suppose that $(x_{\alpha})$ is a net in 
$\Phi_A$ such that every hull-kernel cluster point of $(x_{\alpha})$ is a 
hull-kernel limit point, and such that that there is an $x\in \Phi_A$ to which 
$(x_{\alpha})$ converges in the Gelfand topology. Then the set of hull-kernel 
limit points of $(x_{\alpha})$ is Gelfand connected.} 
\bigskip 
\noindent {\bf Proof.} Let $L$ be the set of hull-kernel limit points of 
$(x_{\alpha})$. Note that $L$ is hull-kernel closed, and hence also
Gelfand closed. Suppose, for a contradiction, that $L$ is a disjoint union
of two non-empty Gelfand closed sets, $M$ 
and $N$ say. Let $S$ and $T$ be disjoint Gelfand open subsets of $\Phi_A$ 
containing $M$ and $N$ respectively. For each $\alpha$, let $K_{\alpha}$ be 
the hull-kernel closure of the set $\{ x_{\beta}: \beta\ge\alpha\}$. Then 
$L=\bigcap_{\alpha} K_{\alpha}$ [\the\Kel; Theorem 2.7], and each $K_{\alpha}$ is Gelfand closed, so a simple 
topological argument shows that there is an $\alpha_0$ such that
for all $\alpha \geq \alpha_0$, $K_{\alpha}\subseteq S\cup T$. If 
we suppose that $x\in M\subseteq S$ then there is a $\gamma \geq \alpha_0$
such that for all $\alpha \geq \gamma$, $x_{\alpha} \in S$.
The quotient Banach algebra $A/I(K_{\gamma})$ has the disconnected maximal 
ideal space $K_{\gamma}$. By the Shilov Idempotent Theorem there exists 
$f\in A$ such that 
$f$ is zero on $K_{\gamma}\cap S$, but $f$ equals one on $K_{\gamma}\cap T$. 
But then $f(x_{\alpha})=0$ for all $\alpha\ge\gamma$. Since the zero set of $f$ is hull-kernel closed, this contradicts the 
hull-kernel-convergence of $(x_{\alpha})$ to points in $N$. Hence 
$L$ is connected. Q.E.D. 
\bigskip 
\noindent {\bf Theorem 3.2} {\sl Let $A$ be a natural, unital Banach function algebra. Let $x\in \Phi_A$. Then 

(i) the hull of $J_x$ is Gelfand connected, 

(ii) the set $\{ y\in \Phi_A: x\sim y\}$ is Gelfand connected.} 
\bigskip 
\noindent{\bf Proof.} (i) Let $E$ be the hull of $J_x$.
Suppose that $y\in E$ with $y\ne x$. 
Then every hull-kernel neighbourhood of $y$ has non-empty intersection with 
every Gelfand neighbourhood of $x$. Thus there is a net $(x_{\alpha})$ in $\Phi_A$ 
converging to $x$ in the Gelfand topology, and to $y$ in the hull-kernel 
topology. By passing to a universal subnet we may assume that every 
hull-kernel cluster point of $(x_{\alpha})$ is a limit point. Thus the set $L$ 
of hull-kernel limit points of $(x_{\alpha})$ is Gelfand connected, by 
Proposition 3.1. 
But if $z\in L$ then $z$ cannot have a hull-kernel neighbourhood disjoint from 
a Gelfand neighbourhood of $x$, so we must have $M_z\supseteq J_x$. 
Thus $L$ is a Gelfand connected subset of $E$, and $y \in L$.
Hence every point of $E$ is in the 
same Gelfand connected component of $E$ as $x$, so $E$ is Gelfand connected. 

(ii) Let $H$ be the set $\{ z\in \Phi_A: x\sim z\}$.
Let $y\in \Phi_A\backslash\{ x\}$ with $y\sim x$. Let $(x_{\alpha})$ be a net 
in $\Phi_A$ converging to both $y$ and $x$ in the hull-kernel topology. By passing 
to a universal subnet we may suppose that $(x_{\alpha})$ is Gelfand 
convergent, and that every hull-kernel cluster point of $(x_{\alpha})$ is a 
limit point. Thus the set $L$ of hull-kernel limit points of $(x_{\alpha})$ is 
connected, by Proposition 3.1. But for each $z\in L$, $z\sim x$. Thus $L$ is a 
Gelfand connected subset of $H$, and $y\in L$.
As above it follows that $H$ is Gelfand connected. Q.E.D. 
\bigskip 
\noindent Theorem 3.2(i) answers one of the questions posed by Gorkin and 
Mortini [\the\GM]. 
\smallskip
One immediate consequence of Theorem 3.2 is the following. 
\bigskip 
\noindent{\bf Corollary 3.3} {\sl Let $A$ be a unital Banach function algebra. Then 
$\Phi_A$ has no isolated points of discontinuity.} 
\bigskip 
We now give examples to show that for a non-natural Banach function 
algebra the hull of $J_x$ does not have to be connected, and there may be isolated
points of discontinuity. 
\bigskip 
\noindent{\bf Examples} Let $A$ be the disc algebra on the countable space $X$ described in the Example after Theorem 2.4. Then the hull of $J_0$ is neither connected, nor uncountable. The only point of continuity is $0$, so every other point of $X$ is an isolated point of discontinuity.

\smallskip
It is possible for a uniform algebra to have a solitary point of discontinuity. 
For example, let $A$ be 
the uniform algebra obtained by restricting $H^{\infty}$ 
(the algebra 
of bounded functions on the disc which are analytic on the open disc) to 
the fibre of its maximal ideal space associated with a point on the unit 
circle, see [\the\Hoff; p.187ff]. Then $A$ is regular on its Shilov boundary 
$\Gamma_A$, but $A$ is not normal. Now consider $A$ as a uniform
algebra on $X=\Gamma_A\cup\{y\}$, where $y$ is any point of 
$\Phi_A\backslash\Gamma_A$. Then $y$ is the solitary point of discontinuity 
for $A$ on $X$ (although there must be many non-R-points --- see the results
in Section 4).
\smallskip
This example also shows that it is possible for the hull of $J_x$ to have 
exactly two points, because if $x$ is any element of $F_y\backslash\{y\}$
then the hull of $J_x$ is precisely the set $\{x,y\}$. An easy modification
produces an example of a uniform algebra on its Shilov boundary with a $J_x$ having a two-point hull. Simply form a new compact space $Y$ by gluing an interval to
$X$ with endpoints at $x$ and $y$ above, and take the uniform
algebra of all continuous functions on $Y$ whose restriction to $X$ is in $A$.
\bigskip 
\noindent Here is another consequence of Theorem 3.2. 
\bigskip 
\noindent{\bf Corollary 3.4} {\sl Let $A$ be a Banach function algebra with 
$\Phi_A= [0,1]$. If $A$ is weakly regular then $A$ is normal.} 
\bigskip 
\noindent{\bf Proof.} By Theorem 2.3, the set of points of continuity contains 
a dense $G_{\delta}$ of $[0,1]$ in the Gelfand topology. Suppose that $x\in 
[0,1]$ with $x$ not an R-point. Then the hull of $J_x$ is a connected subset 
of $[0,1]$ by Theorem 3.2(i), and hence is an interval. Thus it contains 
points of continuity other than $x$ itself, contradicting the definition of a 
point of continuity. Thus every point of $[0,1]$ is an R-point, so $A$ is 
regular, and hence normal. Q.E.D. 
\bigskip 
\noindent If the condition on the character space is dropped, there are non-regular, 
weakly regular uniform algebras on $[0,1]$. For example, let $A$ be the 
non-trivial uniform algebra on the 
Cantor set described in [\the\We; 9.3]. The character space 
of $A$ is the whole of the Riemann sphere (see [\the\We; 9.2, 9.3]), 
so [\the\St; 27.3] shows that $A$ is not normal on the Cantor set. 
Let $B$ be the algebra of continuous functions on $[0,1]$ whose 
restrictions to the Cantor set lie in $A$. Then $B$ is a non-normal 
uniform algebra on $[0,1]$. Each point of $[0,1]$ not in the Cantor set is a 
separated point, so $B$ is weakly regular on $[0,1]$. The character space of 
$B$ is not equal to $[0,1]$ since it contains a copy of the Riemann sphere. 
\bigskip
\noindent Now let $A$ be a Banach function algebra on a compact Hausdorff space $X$, and let $x\in X$. Then for $y$ in the hull of $J_x$ there is a net converging to both $x$ and $y$ in the hull-kernel topology. It is not necessarily the case, however, that a net can be found converging to an arbitrary pair of points in the hull of $J_x$. Consider the following example.
\bigskip
\noindent {\bf Example} Let $A$ be the disc algebra on the disc $X$. Let $x$ be a fixed point in $X$, and let $B=\{ (f,g)\in A\oplus A: f(x)=g(x)\}$. Then the character space of $B$ consists of two copies of the disc glued at the point $x$. The ideal $J_x$ is the zero ideal, so the hull of $J_x$ is the whole of $\Phi_B$. But if $y$ and $z$ belong to different copies of the disc then there is no net in $\Phi_B$ converging simultaneously to both $y$ and $z$ in the hull-kernel topology.
\bigskip
\noindent To show that there are $x$'s for which the hull of $J_x$ is contained in a hull-kernel limit set, we require the notion of a primal ideal. Recall that an ideal $I$ in a commutative ring $R$ is {\sl primal} 
if whenever $a_1, a_2,\ldots,a_n\in R$ with $a_1a_2\ldots a_n=0$, then $a_i\in 
I$ for at least one $i\in\{ 1,\ldots, n\}$, see [\the\MinP] for example. It is a 
straightforward piece of general topology to show that if 
$A$ is a Banach function algebra, and $I$ is a closed primal ideal of $A$ then there is a net in $\Phi_A$ converging to every point in the hull of $I$ in the hull-kernel topology, see [\the\AB; 3.2].

We shall show that if $A$ is a separable Banach function algebra then, in the sense of Baire category, most of the ideals $\overline {J_x}$ are primal. First we observe that every closed primal ideal contains a $J_x$.
\bigskip 
\noindent{\bf Lemma 3.5} {\sl Let $A$ be a Banach function algebra on a 
compact Hausdorff space $X$. Let $P$ be a closed primal ideal of $A$. Then 
there exists 
$x\in X$ such that $P\supseteq J_x$.} 
\bigskip 
\noindent {\bf Proof.} Suppose for a contradiction that $P\not\supseteq J_x$ 
for each $x\in X$. Thus for each $x\in X$ there exists $f\in A$ such that $f$ 
vanishes in a neighbourhood of $x$, but $f\notin P$. By a compactness 
argument, there exist a finite number of functions $f_1,\ldots,f_n$ such that 
$f_1\ldots f_n=0$, with $f_i\notin P$, $(1\le i\le n)$. This contradicts the 
assumption that $P$ is primal. Hence $P$ contains $J_x$ for some $x\in X$. 
Q.E.D. 
\bigskip 
\noindent {\bf Lemma 3.6} {\sl Let $A$ be a Banach function algebra on a 
compact Hausdorff space $X$. For $f\in A$ the function $x\mapsto \Vert 
f+\overline{J_x}\Vert$ is upper semi-continuous on $X$.} 
\bigskip 
\noindent {\bf Proof.} Let $f\in A$, $x\in X$, and let $\epsilon>0$ be given. 
Then there exists $g\in J_x$ such that $\Vert f-g\Vert<\Vert 
f+\overline{J_x}\Vert +\epsilon$. But $g\in J_y$ for all $y$ in a 
neighbourhood $N$ of $x$, so $\Vert f+\overline{J_y}\Vert\le\Vert 
f-g\Vert<\Vert f+\overline{J_x}\Vert+\epsilon$ for all $y\in N$. Thus the norm 
function $x\mapsto\Vert f+\overline{J_x}\Vert$ is upper semi-continuous on $X$. 
Q.E.D. 
\bigskip 
\noindent For an upper semi-continuous function $f$ on a Baire space $X$, the 
set of points of continuity of $f$ contains a dense $G_{\delta}$ of $X$, see 
[\the\Dix; B18] for example. Suppose now that $A$ is a separable Banach 
function algebra on a compact Hausdorff space $X$, and let $\{ 
f_i\}_{i=1}^{\infty}$ be a countable dense subset of $A$. Then $X$ is a Baire 
space, so there is a dense $G_{\delta}$ $Y$ of $X$ consisting of points at 
which all the norm-functions $x\mapsto \Vert f_i+\overline{J_x}\Vert$ are 
continuous. But is is straightforward to check that the continuity of these 
norm-functions, at a particular point, for a dense subset of $A$ forces the 
continuity of all the norm functions at that point. Thus if $A$ is separable 
there is a dense $G_{\delta}$ of $X$ consisting of points at which every norm 
function $x\mapsto \Vert f+\overline{J_x}\Vert$ is continuous. 
\bigskip 
\noindent{\bf Proposition 3.7} {\sl Let $A$ be a separable Banach function 
algebra on a compact Hausdorff space $X$. Then the set of $x\in X$ for which 
$\overline{J_x}$ is primal is a dense $G_{\delta}$ in $X$.} 
\bigskip 
\noindent{\bf Proof.} By the remarks above, there is a dense $G_{\delta}$ 
$Y\subseteq X$ such that for each $y\in Y$ every norm function is continuous 
at $y$. Let $y\in Y$. We show that $\overline{J_y}$ is primal. Suppose that 
$f_1,\ldots ,f_n$ are a finite number of elements of $A$ which are not in 
$\overline{J_y}$. By continuity of the norm functions at $y$, there is a 
neighbourhood $N$ of $y$ in $X$, such that 
$f_i\notin \overline{J_x}$ for $x\in N$. For each $i\in\{ 1,\ldots ,n\}$, let 
$N_i=\{ x\in N:f_i(x)\ne 0\}$. Then $N_i$ is a dense open subset of $N$, so 
$\bigcap_{i=1}^n N_i$ is non-empty. Let $z$ belong to this intersection. Then 
$f_i(z)\ne 0$ for $1\le i\le n$, so $f_1f_2\ldots f_n\ne 0$. This shows that 
$\overline{J_y}$ is primal. Q.E.D. 
\bigskip 
\bigskip
\noindent{\bf 4. Points of discontinuity} 
\bigskip 
\noindent In this section we consider what happens when there is a point $y$ 
of discontinuity. The general aim is to show that $F_y$ has to be big, and 
hence that there have to be many non-R-points. The main results are that if 
$A$ is a natural unital Banach function algebra and $y\notin \Gamma_A$ then $F_y$ has 
an infinite, connected subset, while if $A$ is any uniform algebra then $F_y$ has 
a perfect subset. Thus in both cases there are uncountably many 
non-R-points. A simple example was given after Corollary 3.3 of a non-natural 
Banach function algebra with a solitary non-R-point, but we are unable to say 
whether this phenomenon can occur for a natural Banach function algebra,
or for a Banach function algebra on its Shilov boundary. 
\bigskip 
\noindent {\bf Definition} Let $A$ be a Banach function algebra on a compact, 
Hausdorff space $X$. Then $A$ is {\sl local on X} if the following condition 
holds: a function on $X$ belongs to $A$ if it agrees in a neighbourhood of 
each point of $X$ with an element of $A$. The algebra $A$ is {\sl $2$-local on 
X} if the following 
condition holds: a function on $X$ belongs to $A$ if there are elements $g_1$ and $g_2$ in $A$ so 
that every point of $X$ has neighbourhood on which $f$ agrees with 
either $g_1$ or $g_2$. If $X$ is the character space of $A$ above, then $A$ is said to 
be {\sl local} or {\sl 2-local} respectively. 
\bigskip 
Every normal Banach function algebra is local, and every local Banach function 
algebra is, of course, $2$-local. Most commonly-met uniform algebras are local. For instance if $X$ is a 
compact subset of the plane then the uniform algebras $A(X)$ and $R(X)$ are 
local. All approximately normal
uniform algebras (hence all uniform algebras with character space equal to 
$[0,1]$) are 2-local [\the\Wi; 2.4, 3.1]. It follows 
from [\the\FS; Theorem 11] that the algebra $B$ after Corollary 3.4 is not
$2$-local on $[0,1]$. 
For any compact plane set $X$, the Banach function algebras on $X$ introduced by Dales and Davie in [\the\DD] are local on $X$.
\bigskip
In the first result of this section we do not require a norm on our algebra. In fact the
result is valid for multiplicative sub-semigroups of $C(X)$ (where $C(X)$ denotes the algebra of continuous, complex-valued functions on a compact Hausdorff space $X$).
\smallskip
\noindent {\bf Proposition 4.1} {\sl 
Let $X$ be a compact Hausdorff space and let $A$ be a subalgebra of $C(X)$. Suppose that $A$ is 2-local on $X$. Let $y \in X$, and recall that 
$F_y=\{x\in X: J_x \subseteq M_y\}$. 
Then $F_y$ is connected.} 
\bigskip 
\noindent{\bf Proof.} Suppose, for contradiction, that $F_y$ is not connected. 
Since $F_y$ is closed, we may write $F_y=E_1\cup E_2$ where $E_1$ and $E_2$ 
are non-empty, disjoint, closed subsets of $F_y$. Suppose that $y\in E_2$. 
Choose a 
closed neighbourhood $N$ of $E_1$ such that 
$N \cap F_y = E_1$, and let $K$ be the boundary of $N$. Then $K$ is a compact 
subset of $X\backslash F_y$. Thus, for each 
$t \in K$ we can find an $f \in J_t$ with $f(y)=1$. By compactness, 
multiplying finitely many of these functions together, 
we can find a function $g \in A$ which vanishes on a neighbourhood of $K$ and 
such that $g(y) = 1$ (if $K$ is empty, simply take $g$ to be zero on $N$, and 
$1$ elsewhere). Define $h(t)$ to be $0$ for $t \in N$, and set 
$h(t)=g(t)$ for other values of $t$. Then $h\in A$ since $A$ is 2-local. But 
then $h \in J_x$ and $h(y) = 1$, contradicting the choice 
of $x$. The result follows. Q.E.D. 
\bigskip 

\noindent{\bf Corollary 4.2} {\sl Let $A$ be a 2-local Banach function algebra 
on $[0,1]$. If the hull-kernel topology on $[0,1]$ is first countable then $A$ 
is regular on $[0,1]$. If $\Phi_A= [0,1]$ and the set of points at which the 
hull-kernel topology is first countable is dense in $[0,1]$, then $A$ is 
normal.} 
\bigskip 
\noindent{\bf Proof.} Let $y$ be a point of first countability for the 
hull-kernel topology. Then $F_y$ is a connected set with no interior, by 
Proposition 4.1 and Theorem 2.4, 
and hence is a singleton. Thus $F_y=\{ y\}$, so $y$ is a point of continuity. 
Thus if the hull-kernel topology is first countable, every point is a point of 
continuity, so $A$ is regular on $[0,1]$. 
The proof of the second statement is as in Corollary 3.4. Q.E.D. 
\bigskip 
\noindent Corollary 4.2 shows that non-regular Dales-Davie algebras on $[0,1]$ 
must have hull-kernel topology which is not first countable. 
\bigskip 
\noindent We now use Proposition 4.1 to show that if 
$A$ is a natural unital Banach function algebra and $y\notin \Gamma_A$ then $F_y$ has 
an infinite, connected subset.

Let $B$ be a uniform algebra on its character space $\Phi_B$. Then 
there is a smallest local uniform algebra $C$ on $\Phi_B$ containing $B$, 
which is obtained from $B$ as follows. Set $L^0(B)=B$, and for each ordinal 
$\sigma$ define $L^{\sigma+1}(B)$ inductively to be the uniform closure of the 
functions belonging locally to $L^{\sigma}(B)$. If $\sigma$ is a limit 
ordinal, let $L^{\sigma}$ be the uniform closure of 
$$\bigcup\{L^{\sigma'}:0\le\sigma'<\sigma\}.$$ 
The process must terminate (with $C(\Phi_B)$, if not before), so eventually 
$L^{\sigma}=L^{\sigma+1}$. For this ordinal $\sigma$, set $C=L^{\sigma}$. Then 
evidently $C$ is the smallest local uniform algebra on $\Phi_B$ containing 
$B$. Examples are given in [\the\Sid] showing that, with $\Phi_B$ metrizable,
the process can terminate at any given stage before the first uncountable 
ordinal. 
\bigskip 
\noindent{\bf Lemma 4.3} {\sl Let $\{ A_{\lambda}\}_{\lambda\in\Lambda}$ be an 
increasing family of uniform algebras, all with the same character space $X$. 
Let $B$ be the uniform closure of their union in $C(X)$. Suppose that 
$\Gamma_{A_{\lambda}}=Y$ for all $\lambda$, where $Y$ is a Gelfand closed 
subset of $X$. Then $\Phi_B=X$ and $\Gamma_B=Y$}. 
\bigskip 
\noindent{\bf Proof.} It is easy to see that $\Phi_B=X$, and it is also clear that 
$\Gamma_B\supseteq Y$, since $\Gamma_B$ is a boundary for each $A_{\lambda}$. 
Let $x\in X$ be a strong boundary point for $B$. Let $N$ be any Gelfand 
neighbourhood of $x$. Then there exists $f\in B$ such that $|f(x)|>3/4$ while 
$|f|_{X\backslash N}<1/4$ [\the\St; 7.18]. Hence for some $\lambda$ there 
exists $g\in A_{\lambda}$ such that $|g(x)|>1/2$ while $|g|_{X\backslash 
N}<1/2$. This shows that $Y=\Gamma_{A_{\lambda}}$ meets $N$, and hence that 
$x\in Y$ since $N$ was an arbitrary neighbourhood of $x$. Thus 
$\Gamma_B\subseteq Y$, since the set of strong boundary points is dense in 
$\Gamma_B$ [\the\St; 7.24]. Hence $\Gamma_B=Y$. Q.E.D. 

\bigskip 
\noindent{\bf Lemma 4.4} {\sl Let $A$ be a unital Banach function algebra and let $C$  be the smallest local uniform subalgebra of $C(\Phi_A)$ containing $A$. Then $\Phi_A=\Phi_C$ and $\Gamma_A=\Gamma_C$.} 
\bigskip 
\noindent{\bf Proof.} First, let $B$ be the uniform closure of $A$ in 
$C(\Phi_A)$. It is straighforward that $\Phi_A=\Phi_B$. We now argue as in Lemma 4.3. Since $\Gamma_B$ is a 
boundary for $A$, we have $\Gamma_B\supseteq \Gamma_A$. On the other hand 
suppose that $x\in\Phi_A$ is a strong boundary point for $B$. Let $N$ be any 
Gelfand neighbourhood of $x$. Then there exists $f\in B$ such that 
$|f(x)|>3/4$ while $|f|_{\Phi_A\backslash N}<1/4$ [\the\St; 7.18]. Hence there 
exists $g\in A$ such that $|g(x)|>1/2$ while $|g|_{\Phi_A\backslash N}<1/2$. 
This shows that $\Gamma_A$ meets $N$, and hence that $x\in\Gamma_A$ since $N$ 
was an arbitrary neighbourhood of $x$. Thus $\Gamma_A=\Gamma_B$. 

Now let $C$ be the smallest local uniform subalgebra of $C(\Phi_A)$ containing 
$B$. Then Stolzenberg showed that $\Phi_{L^{\sigma}(B)}=\Phi_{L^{\sigma+1}(B)}$ for each ordinal $\sigma$ [\the\Stolz], and he also mentioned in a remark that $\Gamma_{L^{\sigma}(B)}=\Gamma_{L^{\sigma+1}(B)}$ for each ordinal $\sigma$
(this follows, as with the maximal ideal space space result, by considering the intermediate uniform algebras generated
by $L^{\sigma}(B)$ together with a single function $f$ which is locally in $L^{\sigma}(B)$).
Thus it follows from Lemma 4.3 that $\Phi_C=\Phi_B$ and $\Gamma_C=\Gamma_B$, and hence that $\Phi_C=\Phi_A$ and $\Gamma_C=\Gamma_A$. Q.E.D. 
\bigskip 
\noindent We are now ready for the main result of this section. 
\smallskip 
\noindent {\bf Theorem 4.5} {\sl Let $A$ be a unital Banach function algebra and let 
$y\in \Phi_A\backslash\Gamma_A$. Then $F_y$ has an infinite, connected 
subset.} 
\bigskip 
\noindent{\bf Proof.} Let $C$ be the smallest local uniform subalgebra of 
$C(\Phi_A)$ containing $B$. Then we saw in Lemma 4.4 that $\Phi_C=\Phi_A$ and 
that $\Gamma_C=\Gamma_A$. Hence $y\notin \Gamma_C$, so Lemma 2.1(ii) shows 
that $y$ is not a point of continuity for $C$. Thus the set $F_{y,C}=\{x\in 
\Phi_C: J_x \subseteq M_y\}$ (where $J_x$ and $M_y$ are here defined relative 
to the algebra $C$) has more than one point, and is connected by Proposition 
4.1. Since $F_y\supseteq F_{y,C}$, the result follows. Q.E.D. 
\bigskip 
\noindent For uniform algebras a weaker but more general statement is true. 
\smallskip 
\noindent {\bf Theorem 4.6} {\sl Let $A$ be a uniform algebra on a compact 
Hausdorff space $X$ and let $y \in X$. Then either $y$ is a point of continuity or $F_y$ 
has a non-empty, perfect subset.} 
\bigskip 
\noindent {\bf Proof.} Suppose, for a contradiction, that $F_y$ has no 
non-empty, perfect subsets, and that there exists $x \in F_y\backslash\{y\}$. 
Choose 
$f \in A$ with $f(x) = 0$ and $f(y) = 1$. By [\the\Ru], $f(F_y)$ is a 
countable, 
compact set. Thus it is easy to choose two disjoint, closed rectangles 
$R_1$ and $R_2$ in the complex plane such that $0\in R_1$, $1 \in R_2$, 
and such that 
$f(F_y)$ is contained in the interior, $U$, of $R_1 \cup R_2$. By Runge's 
theorem [\the\Rud; 13.9] we can choose a sequence $(p_n)$
of polynomials such that $p_n \to 0$ uniformly on $R_1$, but $p_n 
\to 1$ 
uniformly on $R_2$, as $n\to\infty$. 
Also, $X\backslash f^{-1}(U)$ is a compact subset of $X\backslash F_y$, so 
by 
compactness (as in Proposition 4.1) 
we can find a function $g$ 
in $A$ which vanishes on a neighbourhood of $X\backslash f^{-1}(U)$, but 
such 
that $g(y) = 1$. 
Set $h_n = g p_n(f)$. Then $h_n \in A$, and the sequence $(h_n)$ converges 
uniformly on 
$X$ to 
a function which vanishes on a neighbourhood of $x$ but which is $1$ at 
$y$. This 
contradicts the choice of $x$. Thus $F_y\backslash\{y\}$ must be empty, as 
claimed. Q.E.D. 
\bigskip 
\noindent The non-trivial uniform algebra $A$ on the Cantor set, mentioned after Corollary 3.4, is not regular on the Cantor set, because it is an 
integral domain. Thus for a point $y$ of discontinuity, the set $F_y$ is 
not connected. 
This example is, however, 2-local on its character 
space (which is in fact equal to the Riemann sphere) so the set $F_y$ 
would be connected if we were working on the character space. 
\bigskip 
\noindent The next corollary is used in [\the\FS]. 
\smallskip 
\noindent{\bf Corollary 4.7} {\sl 
Let $A$ be a uniform algebra on a compact Hausdorff space $X$. Suppose that $A$ is not regular on 
$X$. Then the set of non-R-points has a non-empty perfect subset.} 
\bigskip 
\noindent{\bf Proof.} Since $A$ is not regular on $X$, there exists $y\in X$ 
such that the set 
$\{x\in X\backslash \{y\}: J_x \subseteq M_y\}$ 
is non-empty. But all points of this set are non-R-points, and, by 
Proposition 4.6, this set contains a non-empty, perfect set. The 
result follows. Q.E.D. 
\bigskip 
\noindent We saw in Section 3 that a Banach function algebra can have a 
solitary non-R-point. We do not know, however, whether this can happen for a 
natural Banach function algebra, nor for a Banach function algebra on its 
Shilov boundary. Proposition 4.1 shows that such an example would have to be 
non-2-local. The next theorem shows that the difficult case is when every 
point is an independent point. 
\bigskip 
\noindent{\bf Theorem 4.8} {\sl Let $A$ be a Banach function algebra on a compact Hausdorff space $X$ and let 
$y$ be a non-independent point of $X$. Then $F_y$ has a non-empty, perfect 
subset.} 
\bigskip 
\noindent{\bf Proof.} Let $B$ be the uniform closure of $A$ in $C(X)$. Then it 
is easy to see that $y$ is not a point of continuity for $B$ on $X$, so the 
set $F_{y,B}=\{x\in X: J_x \subseteq M_y\}$ (where $J_x$ and $M_y$ are here 
defined relative to the algebra $B$) has an infinite perfect subset by Theorem 
4.6. Since $F_y\supseteq F_{y,B}$, the result follows. Q.E.D. 

\bigskip 
\noindent Recall that a Banach function algebra $A$ on a compact Hausdorff 
space $X$ is {\sl approximately regular on $X$} if whenever $Y$ is a closed 
subset of $X$, and $x\in X\backslash Y$, then for any $\epsilon>0$ there 
exists $f\in A$ such that $f(x)=1$ and $|f|_Y<\epsilon$. This is clearly 
equivalent to every point of $X$ being an independent point. It is also clear 
that if $A$ is a uniform algebra and $A$ is approximately regular on $X$ then 
$X$ must be the Shilov boundary of $A$. 

The disc algebra is approximately regular on the circle, but the tomato-can 
algebra is not approximately regular on its character space. 

\bigskip 
\noindent {\bf Lemma 4.9} {\sl Let $A$ be a uniform algebra on a compact, 
Hausdorff space $X$. Let $x\in X$ be an independent point, and let $y\in X$. 
Then $x$ cannot be an isolated point in $F_y$, unless $x=y$ and $y$ is a point 
of continuity.} 
\bigskip 
\noindent{\bf Proof.} Suppose, for a contradiction, that $x$ is an isolated 
point of $F_y$ and that $F_y\backslash\{x\}$ is non-empty. Since $x$ is an 
independent point there is a function $f'\in A$ taking the value $1$ at $x$ 
and nearly zero on $F_y\backslash\{x\}$. Hence $f=1-f'$ is zero at $x$ and 
nearly $1$ on $F_y\backslash\{x\}$. We now argue as in Theorem 4.6. 

It is easy to choose two disjoint, closed rectangles 
$R_1$ and $R_2$ in the complex plane such that $0\in R_1$, $1 \in R_2$, 
and such that 
$f(F_y)$ is contained in the interior, $U$, of $R_1 \cup R_2$. By Runge's 
theorem [\the\Rud; 13.9] we can choose a sequence $(p_n)$
of polynomials such that $p_n \to 0$ uniformly on $R_1$, but $p_n 
\to 1$ 
uniformly on $R_2$. 
Also, $X\backslash f^{-1}(U)$ is a compact subset of $X\backslash F_y$, so 
by 
compactness (as in Proposition 4.1) 
we can find a function $g$ 
in $A$ which is $0$ on a neighbourhood of $X\backslash f^{-1}(U)$, but 
such 
that $g(y) = 1$. 
Set $h_n = g p_n(f)$. Then $h_n \in A$, and the sequence $(h_n)$ converges 
uniformly on 
$X$ to 
a function $h$ which is $0$ on a neighbourhood of $x$ but which is $1$ on a 
neighbourhood of $F_y\backslash\{x\}$. 

The function $h$ shows that $x=y$, for otherwise we have a function vanishing 
in a neighbourhood of $x$, but non-zero at $y$, which
contradicts the fact that $x \in F_y$. Q.E.D. 
\bigskip 
\noindent As an immediate consequence of Lemma 4.9, we have the following. 
\smallskip 
\noindent {\bf Theorem 4.10} {\sl Let $A$ be a uniform algebra. If $A$ is 
approximately regular on $X=\Gamma_A$, and $y\in \Gamma_A$, then $F_y\backslash\{y\}$ has no 
isolated points. Hence $A$ has no isolated non-R-points.} 
\bigskip
We conclude with an example showing that a uniform algebra on its
Shilov boundary can have an isolated non-R-point.
We do not know if this is possible for a uniform algebra on its
character space.
\bigskip 
\noindent
{\bf Example} Let $X$ be the closed unit disc and let $Y$ be the
set $T\cup\{0, {1\over 2}, {1\over 3},\dots\}$. Let $A$ be the uniform algebra of
all functions in $C(X)$ whose restriction to $Y$ is in
the restriction to $Y$ of the disc algebra. Then it is easy to see
that $X$ is the Shilov boundary of $A$, and that the only non-R-points
for $A$ are the points of $T$ and the point $0$. Thus $0$ is an isolated
non-R-point for $A$. In fact, for $y \in Y$, $F_y = \{0,y\}\cup T$.
All other points of $X$ are points of continuity for $A$.
\bigskip 
\bigskip
\centerline{\bf References} 
\medskip 
\input TempReferences.tex 
\bigskip 
\bigskip

\centerline{School of Mathematical Sciences} 
\centerline{University of Nottingham} 
\centerline{ NG7 2RD} 
\centerline{ U.K.} 
\medskip \centerline{email: Joel.Feinstein@nottingham.ac.uk} 
\bigskip 
\centerline{Department of Mathematical Sciences} 
\centerline{University of Aberdeen} 
\centerline{AB24 3UE} 
\centerline{U.K.} 
\medskip 
\centerline{e-mail: ds@maths.abdn.ac.uk}

\end